\documentclass[10pt]{amsart}
\usepackage{amsopn}
\usepackage{amsmath, amsthm, amssymb, amscd, amsfonts, xcolor}
\usepackage[latin1]{inputenc}
\usepackage{xcolor}
\usepackage{enumitem}
\usepackage{multirow}

\input xy
\xyoption{all}
\usepackage[OT2,OT1]{fontenc}
\newcommand\cyr{%
\renewcommand\rmdefault{wncyr}%
\renewcommand\sfdefault{wncyss}%
\renewcommand\encodingdefault{OT2}%
\normalfont
\selectfont}
\DeclareTextFontCommand{\textcyr}{\cyr}

\usepackage[pdftex]{hyperref}

\theoremstyle{plain}

\newtheorem{teo}{Theorem}[section]
\newtheorem{cor}[teo]{Corollary}

\newtheorem{prop}[teo]{Proposition}
\newtheorem{lema}[teo]{Lemma}

\theoremstyle{definition}
\newtheorem{defi}[teo]{Definition}

\newtheorem{nota}[teo]{Remark}

\numberwithin{equation}{section}

\makeatletter\renewcommand\theenumi{\@roman\c@enumi}\makeatother

\title{Normal holonomy and rational properties of the  shape operator}

\author{Carlos Olmos}
\address{Facultad de Matem\'atica, Astronom\'ia y F\'isica, Universidad Nacional de C\'ordoba, 
Ciudad Universitaria, 5000 C\'ordoba, Argentina.}
\email{olmos@famaf.unc.edu.ar}

\author{Richar Ria\~ no-Ria\~ no }
\address{Departamento de Matem\' aticas, Facultad de Ciencias, Universidad De Los Andes
	Cra 1\#18A-12 Bogot\' a, Colombia.	}
\email{rf.riano@uniandes.edu.co}
\thanks{Supported by Famaf-UNC, CIEM-Conicet, Argentina/ Colciencias and Universidad de Los Andes, Colombia}
\thanks {This research started   during the visit of the first author to the  Universidad de Los Andes  and was essentially finished during the visit of the second author to FaMAF}


\begin {document}

\begin{abstract} Let $M$ be a most singular orbit of the isotropy representation of a simple symmetric space.  Let $(\nu _i, \Phi _i)$ be an irreducible factor of the normal holonomy representation 
	$(\nu _pM, \Phi (p))$. 
	 We prove that there exists a basis of a  section $\Sigma _i\subset \nu _i$ of $\Phi _i$ such that the corresponding shape operators have rational eigenvalues (this is not in general true for other isotropy orbits). Conversely,  this property, if referred to some non-transitive irreducible normal holonomy factor,  characterizes the isotropy orbits. 
	 We also prove that the definition of a submanifold with constant principal curvatures can be given by using only the traceless shape operator, instead of the shape operator, restricted to  a non-transitive (non necessarily irreducible) normal holonomy factor.  
	 This article  generalizes previous results of the authors that characterized Veronese submanifolds in terms of normal holonomy. 
	 
\end{abstract}

\maketitle 

\section {Introduction}

The orbits of $s$-representations (i.e.,  isotropy representations of semisimple  symmetric spaces) play a similar role, in Euclidean submanifold geometry, 
 to that of symmetric spaces in Riemannian geometry\cite {BCO}. 
These orbits coincide with the so-called generalized real flag manifolds. 

Any $s$-representation is a polar representation, i.e., there is a linear subspace, called a section (which in this case is a maximal abelian subspace of the Cartan complement), that intersects every orbit in an orthogonal way. Conversely, by the classification of J. Dadok \cite {Da, EH}, given a polar representation there exists an $s$-representation with the same orbits. 

The rank of  an $s$-representation, i.e.,    the codimension of a principal orbit,  coincides with the rank of the associated symmetric space. 

Orbits of $s$-representations are submanifolds with constant principal curvatures. Namely,  the shape operator of such an orbit, in the direction of  any  arbitrary  parallel normal field  along a curve, has constant eigenvalues. Moreover, the principal orbits are isoparametric submanifolds, i.e., submanifolds with constant principal curvatures and flat normal bundle.
Conversely, by a remarkable result of Thorbergsson \cite {Th, O1}, any 
 irreducible and full Euclidean isoparametric submanifold of codimension at least $3$ is a principal orbit of  an $s$-representation.

The  Thorbergsson theorem was genereralized  by A. J. Di Scala and the  first author \cite {DO} to the so-called  rank rigidity theorem for submanifolds. In particular, one has the following characterization, where the rank of a submanifold is the maximal number of linearly independent, locally defined, parallel normal fields.  

\

\begin {teo}\label {1.1} [\cite {O2, O3}] \label {RR} Let $M ^n$, $n\geq 2$,  
be an irreducible and full homogeneous submanifold of Euclidean space. 
Then 

\ {\rm (i)} $\mathrm {rank}(M)\geq 1$ if and only if $M$ is contained in a sphere.

{\rm (ii)} $\mathrm {rank} (M)\geq 2$ if and only if $M$  is an orbit of an $s$-representation, of rank at least $\mathrm {rank} (M)$,  which is not  most singular.
\end {teo}

\

Since the normal holonomy is a conformal invariant, any conformal diffeomorphism of the sphere maps a submanifold into another one with the same rank. But, in general, such a diffeomorphism does not preserve the homogeneity. So the assumption that the submanifold is homogeneous   cannot be dropped from the above result. The group of presentation of $M$, in the above theorem, may be smaller than the group associated to the $s$-representation.

Observe that the rank of a submanifold can be regarded as the dimension of the fixed vectors of the normal holonomy action. 
So, the rank rigidity theorem characterizes, in terms of the normal holonomy, all  the orbits of $s$-representations which are not most singular. 

 A major open question,  in this context, is the following conjecture that would generalize the rank rigidity theorem. 
 
 \

 \noindent  {\bf Conjecture}(\cite {O2}) \it  A full irreducible homogeneous submanifold $M ^n$ of the sphere,   $n\geq 2$,  such that the normal holonomy is non-transitive must be  an orbit of an $s$-representation.
\rm

\

 A partial  answer to this conjecture was given in \cite {OR}.  Namely, the so-called Veronese submanifolds of dimension at least $3$ are characterized as homogeneous submanifols with generic first normal space and  irreducible  non-transitive normal holonomy. In particular, this implies a positive answer for the  conjecture if $n= 3$, since the number of irreducible factors of the normal holonomy  is always  bounded by $\frac n2$. For $n=2$ the normal holonomy must be always transitive and so the conjecture also holds. 

\

The  goal of this paper is to generalize the results in  \cite {OR}. For this sake we need, in particular, to avoid the delicate topological arguments used there for $n=3$.  This led us, on the one hand,  to consider the so-called traceless shape operator $\tilde A$ and to recover, from this object,   the information given by the usual shape operator.

On the other hand, we  need  necessary and sufficient conditions 
on the subspace $\tilde A _\Sigma\subset \mathrm {End}(T_pM)$, $\Sigma$ a section for the normal holonomy action, so that $M$ is an orbit of an $s$-representation. 
So, we first study how this subspace looks like for most singular $s$-orbits.
We prove that it has interesting rational properties.
Such properties, that depend on the  fact that the Weyl group associated to an $s$-representation is crystallographic,  do  not hold for arbitrary isotropy orbits which are not most singular. Namely, such a subspace admits a basis any of whose elements has rational eigenvalues. The same properties are true if we replace $\Sigma$ by a section of the normal holonomy action restricted to an irreducible factor. 

The family of shape operators at $p$ of a (Euclidean or spherical) submanifold $M^n$ may be regarded as a subspace of $T_{[e]}X$, where $X$ is the symmetric space $ \mathrm {GL}(T_pM)/\mathrm {SO}(T_pM) \simeq \mathrm {GL}_n/\mathrm {SO}_n$. Since this symmetric space is not semisimple, it is natural to consider  the so-called {\it traceless shape operator}
$\tilde A_\xi = A_\xi - \frac 1n\langle H , \xi\rangle Id = A_\xi - \frac 1n\mathrm {trace} (A_\xi) Id$, where $A$ is the shape operator and $H$ is the mean curvature vector.  In this way, $X$ is replaced by the simple symmetric space 
$ \mathrm {SL}(T_pM)/\mathrm {SO}(T_pM)$. Moreover, if we replace the shape operator $A_\xi$ by the so-called {\it dual shape operator} $iA_\xi$, the above mentioned rational properties are equivalent to the fact that a maximal abelian subspace of shape operators is the tangent space to a compact flat, in general non-maximal,  of the dual symmetric space 
$X^* \simeq  \mathrm {SU}_n/\mathrm {SO}_n$. This allows us to use general facts about smooth variations of compact flats. 

Let us enounce our main results

\begin {teo} \label {traceless} Let $M^n$, $n\geq 2$,  be a full and irreducible  submanifold either of the sphere $S^{N -1}$ or  of the Euclidean space  $\mathbb R^N$ which is not contained in a sphere. Let $0\neq \nu 'M$ be a parallel subbundle of the normal bundle $\nu M$ of $M$. 
Assume that the normal holonomy group $\Phi (p)$ of $M$, restricted to $\nu '_pM$, is not transitive on the sphere, $p\in M$ (e.g.,   if $M$  is not a hypersurface and has flat normal bundle, and $\nu'M=\nu M$). Then the traceless shape operator of $M$, restricted to $\nu' M$,  has constant principal curvatures if and only if $M$ is a submanifold with constant principal curvatures. 
\end {teo}

\

Let us say that the traceless shape operator $\tilde A$ of $M$, restricted to $\nu' M$,  has constant principal curvatures if $\tilde A_{\xi (t)}$ has constant eigenvalues for any parallel section, along any curve, $\xi (t)$ of $\nu' M$.

This result  will be needed for the proof of the next theorem. Moreover, it implies  substantial simplifications in the proof of the main result in \cite {OR}, that avoids topological arguments.

The above theorem is not true if we drop the condition that the normal holonomy is not transitive. In fact, one can construct  a Weingarten, non isoparametric, rotational  surface in $\mathbb R ^3$  such that the traceless   shape operator has constant eigenvalues. Or, equivalently, the difference between the two principal curvatures is constant.
 The construction of such an example is just by solving an ordinary differential equation.

The proof of the above theorem, since it is independent of the main tools developed in this article, will be given in the last section.

\

\begin {teo} \label {main}Let $M^n$, $n\geq 2$, be a full and irreducible homogeneous submanifold of the sphere $S^{N-1}$ such that the normal holonomy has  a non-transitive irreducible factor $(\nu _i  , \Phi _i)$. Then the following assertions are equivalent: 
\begin {itemize}
\item [a)] $M$ is an orbit of an irreducible $s$-representation. 
\item [b)] The normal holonomy factor $( \nu _i  ,  \Phi _i)$  has sections of  compact type.
\end {itemize}

\end {teo}

The normal holonomy factor $( \nu _i  ,  \Phi _i)$  has {\it sections of  compact type} if the sections of $\Phi _i$ correspond, via the dual traceless shape operator, to the tangent space of a compact flat in $X^* \simeq \textrm {SU}_n/\textrm {SO}_n$. Or, equivalently, if such a section, via the traceless shape operator $\tilde A$, admits a basis any of whose elements has rational eigenvalues.
If $ \Phi _i $ is non-transitive, on the unit sphere of $ \nu _i $, then  $\Phi _i$ has sections of the compact type if $\tilde A_{\nu _i}$ is a Lie tryple system of the symmetric endomorphisms of the tangent space. This is always the case when $M$ is extrinsic symmetric or when  the normal holonomy acts irreducibly, non-transitively, and the codimension is the maximal one $\frac 12 n(n-1) -1$. Then our result generalizes the main results in \cite {OR} that characterize Veronese submanifolds. 

Any  most singular orbit $M$, of the isotropy representation of 
a simple symmetric space of rank at least $3$,  has a non-transitive irreducible normal holonomy factor. In fact, this follows from  \cite {HO}, and from the table of Dynkin diagrams with root multiplicities for irreducible
Riemannian symmetric spaces.  
Then the above theorem is a characterization of such singular orbits. 

For the same  reasons mentioned after Theorem \ref {1.1}, the assumption of homogeneity cannot be dropped from the hypothesis of Theorem 1.3.
In fact, if one applies a conformal diffeomorphism, 
the traceless shape operator 
changes by a common scalar multiple and the normal holonomy is preserved. 

We hope that this theorem would be useful for proving the previously  mentioned conjecture. 

\section {Preliminaries} \label {Pre}

Let $\mathbb V$ be a $k$-dimensional Euclidean vector space and consider the global symmetric space 
$$X:= \mathrm {SL}(\mathbb V)/\mathrm {SO} (\mathbb V)\simeq  \mathrm {SL}_k/\mathrm {SO}_k, $$ with    Cartan decomposition 
$$ \mathfrak {sl}(\mathbb V)= \mathfrak {sim}^0(\mathbb V) \oplus \mathfrak {so}(\mathbb V),$$ 
where $\mathfrak {sim}^0(\mathbb V)$ are the trace-free symmetric endomorphisms.                   
Let us consider the dual compact globally symmetric space 
   $$X ^*:= \mathrm {SU}(\mathrm V^{\mathbb C})/\mathrm {SO} (\mathrm V)\simeq  \mathrm {SU}_k/\mathrm {SO}_k,$$ with    Cartan decomposition 
   $$\mathfrak {su}(\mathbb V)= i \, \mathfrak {sim}^0(\mathbb V) \oplus \mathfrak {so}(\mathbb V),$$                      where $\mathrm V^{\mathbb C}= \mathbb V \oplus i \mathbb V $ is the complexification of $\mathbb V$. 
There is a natural linear map $\ell$, which maps Lie triple systems into Lie triple systems, from the Cartan subspace $\mathfrak p : = \mathfrak {sim}^0(\mathbb V)$ into the Cartan subspace $i    \mathfrak p = \mathfrak {su}(\mathbb V)  $. Namely, $\ell (w)= iw$. 

 The exponential maps: $\exp : T_{[e]}X \simeq \mathfrak p$, 
 $\exp * : T_{[e]}X^* \simeq i \mathfrak p$ is given by the usual exponential $\mathrm e^x$ of endomorphisms (projected to the respective quotients). 
 
 An abelian subspace $i \mathfrak a$ of $i \mathfrak p\simeq T_{[e]}X^*$ is called of 
 {\it compact type} if it is the tangent space at $[e]$ of a compact  (non-necessarily maximal) flat of $X ^*$. This is equivalent to the fact that 
 $\mathrm e^{i \mathfrak a}$ is compact. 
 
 \
 
 \begin {defi} \label {compact-type} An abelian subspace  $\mathfrak a$ of $\mathfrak p = \mathfrak {sim}^0(\mathbb V)$ is called of {\it compact type} if $\mathrm e ^{i\mathfrak a}$ is compact.
\end {defi}

\begin {nota} \label {curvature}

The Riemannian curvature tensor of $X^*$ at ${[e]}$, up to a rescaling of the metric and identifying $T_{[e]}X^*\simeq \mathfrak p ^*$, is given by 
$$\langle R_{u,v}w, z\rangle =- ( [iu, iv], [iw, iz]), $$
where $(\, , \ )$ is the usual inner product of symmetric matrices given by 
$$(A,B) = \textrm {trace} (AB)$$

\end {nota}

\subsection {Smooth variation of compact flats in symmetric spaces}

\

We will need the following auxiliary  result which is standard to prove. 

\begin {lema} \label {orbit-type}  Let $H$ be a compact connected subgroup of $\text {SO}_r$ and let $c(t)$ be a smooth curve in $\mathbb {R}^r$ such 
that the orbits $H\cdot c(t)$ are all of the same isotropy type, for all $t$. Then, if $p=c(0)$, there exist a smooth curve $\alpha (t)$ in 
the normal space $\nu _{p}(H\cdot p)$ and  a smooth curve $h(t)$ in $H$, with $h(0) = e$, such that:
\begin {itemize}
 \item [i)]$c(t) = h(t)\alpha (t)$, 
 \item [ii)]
$\alpha (t)$ is fixed by the isotropy $H_p$, for all $t$.
\end {itemize}
\end {lema}

Let $ G'/K'$ be a compact 
Riemannian globally   symmetric  space, where $(G',K')$ is an almost effective symmetric pair  with Cartan decomposition 
$$\mathfrak g' = \mathfrak p' \oplus \mathfrak k'.$$ We do not assume that $G'/K'$ is simply connected but we assume that $K'$ is connected. We will be interested, for the applications,  in the symmetric space $\mathrm {SU}_n/\mathrm {SO}_n$. 

The group $K'$ acts on $\mathfrak p'$ by means of the adjoint representation, which is identified with the isotropy representation of $K'$ on $T_{[e]}(G'/K')$.

 Let $\mathfrak {a}'_t \subset i \mathfrak {p}'$, $t \in (-\varepsilon, \varepsilon)$, be an abelian subspace which is tangent to a compact flat of $G'/K'$. Let us assume that this is a  smooth linear variation. Namely, there exists a linear isomorphism $F_t: \mathfrak a' _0 \to \mathfrak a_t $ such that $F_t (w)$ is a smooth curve  in $\mathfrak p '$, for all $w\in \mathfrak a '_0$. 
From the linearity, we  only need to consider $w$ in an open non-empty subset of $\mathfrak a ' _0$. This is  equivalent to the fact that 
$F:(-\varepsilon , \varepsilon)\times \mathfrak a '_0 \to \mathfrak p'$ is smooth. 
Such a smooth variation is called a {\it compact  variation of 
	abelian spaces}.

\begin {prop}\label {mainpro}

 Let $c(t)$ be a smooth curve in $\mathfrak p'$ 
 with $c(0) = p$ and let $\mathfrak a'_t$ be a compact variation of 
abelian subspaces of $\mathfrak p '$ such that 
$c(t)\in \mathfrak a'_t$ and $\dot {c}(t)\perp \mathfrak a '_t$.  Assume, furthermore,  that 
$K'\cdot c(t)$ are all of the same isotropy type and write, as in Lemma \ref {orbit-type}, $c(t)= h(t)\alpha (t)$, where $h(t)$ is a  curve in $K'$ and 
$\alpha (t)$ is a  curve in $\nu _{p}(K'\cdot p)$ 
 fixed by $K'_p$.
Then $\alpha (t) \equiv p$ and hence  the curve $c(t) = h(t)p$ lies in the orbit 
$K'\cdot p$. 
\end {prop}

\begin {proof}  Let $\nu _0$  be the set of vectors fixed by $K'_p$ in 
$\nu _{p}(K'\cdot p)$. One has that $\nu _0\subset \mathfrak p'$ is an abelian subspace of  compact type. In fact, let  $\xi \in \nu _p (K'\cdot p)$ be  such that 
$M= K'\cdot q$ is a principal orbit (and hence isoparametric), where $q=p  + \xi$.  Then, by the  Slice Theorem of Palais-Terng  \cite [Theorem 6.5.9, (i), p. 134] {PT},  
$\nu _0$ coincides with the intersection of 
the normal spaces $\nu _{r}M$, where 
$r\in S_{q,-\xi}: p + (K'_p)^o\cdot \xi$. Since any of these normal spaces is a maximal abelian subspace of $\mathfrak p'$ (and hence of  compact type), we conclude that $\nu _0$ is an abelian subspace of  compact type. 

Let us now consider  the abelian subspace $h^{-1}(t)a'_t$ which is of  compact type and  is contained in $\nu _{\alpha (t)}(K'\cdot \alpha (t))$.  Since 
$K'\cdot p$ and $K'\cdot \alpha (t)$ are parallel manifolds, $\nu _p (K'\cdot p)= \nu _{\alpha (t)}
(K'\cdot \alpha (t))$, for all $t$. 
So 
$$h^{-1}(t)a'_t\subset \nu _p (K'\cdot p).$$
The intersection $h^{-1}(t)a'_t\cap \nu _0$ is an abelian subspace of  compact type that contains $\alpha (t)$. For an open and dense subset $U$ of the real parameter $t$ one has  that $h^{-1}(t)a'_t\cap \nu _0$ has locally constant dimension. Since any abelian subspace of compact type has a discrete number of subspaces of  compact type, the family $h^{-1}(t)a'_t\cap \nu _0$ is locally constant, $t\in U$. 
Let $t_0\in U$ and let $I_{t_0}$ be the open real interval defined by the connected component of 
$t_0$ in $U$. Let 
$$\mathfrak b_{t_0}: = 
h^{-1}(t)a'_t\cap \nu _0,$$
which is independent of $t\in I_{t_0}$. Then, if $t\in I_{t_0}$, 
\begin {equation}\label{(1.1)}
\begin {split}
\dot {c}(t)  &= (h(t)\alpha (t))' = 
\dot {h}(t) \alpha (t) + h(t) \dot {\alpha }(t)
\\ &= h(t) (X_t.\alpha (t)) + h(t) \dot {\alpha }(t),
\end {split}
\end {equation}
where $X_t = h^{-1}(t)\dot {h}(t) \in \mathfrak k'$.
 
On the one hand,  $$X_t.\alpha (t)\in T_{\alpha (t)}
(K'\cdot \alpha (t)) = T_p(K'\cdot p)\perp  
 \nu _{p}(K'\cdot p)$$
On the other hand, the abelian subspace 
$h^{-1}(t)a'_t$ contains $\alpha (t)$ and so it is perpendicular to the orbit $K'\cdot \alpha (t) = 
K'\cdot p$. This implies  that the first term, in  the last equality of \eqref{(1.1)}, is perpendicular to 
$\mathfrak a'_t$. So, from the assumptions of the proposition, the second term of this equality must be perpendicular to $\mathfrak a ' _t$. 

\noindent But  ${\alpha (t)}\in 
\mathfrak b _{t_0}$ and so 
$\dot {\alpha (t)} \in 
\mathfrak b _{t_0}\subset h^{-1}(t)a'_t$. 
Then $h(t) \dot {\alpha (t)}\in a'_t$. A contradiction, unless $\dot {\alpha (t)}=0$. 
for all $t\in U$. Since $U$ is dense, we obtain that $\dot {\alpha (t)}=0$ for all  $t$ and hence 
$\alpha (t) \equiv p$.

\end {proof}

We will be interested in the symmetric spaces of Section \ref {Pre}. Namely, 
$X= \mathrm {SL}(\mathbb V)/\mathrm {SO} (\mathbb V)\simeq  \mathrm {SL}_k/\mathrm {SO}_k, $
and its dual $X ^*= \mathrm {SU}(\mathrm V^{\mathbb C})/\mathrm {SO} (\mathrm V)\simeq  \mathrm {SU}_k/\mathrm {SO}_k$.

We keep the notation of that section. 

\

\begin {defi} A  family $\mathfrak a_t$ of abelian subspaces of $\mathfrak p: = \mathfrak {sim}^0(\mathbb V)$ is said to be a {\it compact variation of abelian subspaces} if 
$ i\mathfrak a_t$ is so. 
\end {defi}

\begin {nota} \label {Prop} From the previous definition,  since the isotropy representations  of dual symmetric spaces coincide, one has that Proposition \ref {mainpro} remains true if we replace $\mathfrak p'$ by $\mathfrak p = \mathfrak {sim}^0(\mathbb V)$. 

\end {nota}

 \

A general reference for the next  part is \cite {BCO}. Let $M ^n =H\cdot v$ be a homogeneous submanifold of the Euclidean space $\mathbb R ^N$, where $H$ is a (connected) Lie subgroup of the isometries $\mathrm {I} (\mathbb {R} ^N)$. Let us denote by $\Phi(u)$ the (restricted) normal holonomy group of $M$ at $u$, which acts on the normal space $\nu _u M$, up to the set of fixed vectors, as an $s$-representation (i.e., the isotropy representation of a semisimple Riemannian symmetric space). In particular, $\Phi (u)$ acts polarly, i.e., there is a linear subspace $\Sigma$ of $\nu _u M$ such that it meets every $\Phi (u)$-orbit and these orbits intersect $\Sigma$ perpendicularly. Such a subspace $\Sigma$ is called a {\it section} of $\Phi (u)$ and is always obtained as the normal space to a principal orbit (and so any other section is conjugated to $\Sigma$ by an element of the normal holonomy group). 

Let us decompose, as in the normal holonomy theorem, 
$$\nu _u M = \nu _0\oplus \cdots \oplus \nu _r\, ,$$
$$\Phi (u) = \Phi _1\times \cdots \times \Phi _r\, ,$$
where $\Phi _i$ acts irreducible on $\nu _i$ and trivially on $\nu _j$, $i=1, \cdots , 
r$, $j= 0, \cdots , r$, $i\neq j$. 
If $\Sigma _i\subset \nu _i$ is a section for the irreducible $s$-representation 
$(\nu _i, \Phi _i)$, then $\Sigma =\Sigma _0 \times \Sigma  _1 \times \cdots \times \Sigma _r$ is a section for the normal holonomy action, where $\Sigma _0 = \nu _0$.  Moreover, any section can be written in this way. 

From Theorem \ref {RR}we may assume that $M = H\cdot v$ is a submanifold of the sphere and that $\nu _0 = 0$.

 If $\xi , \eta \in \Sigma $ then $[A_\xi , A_\eta] =0$, where $A$ denotes the shape operator of $M$. This is a consequence of the fact that the (principal) holonomy tubes have flat normal bundle and the so-called tube formula. Then $A_{\Sigma}$ is an abelian subspace of $\mathfrak {sim}(T_uM)$, the symmetric endomorphisms of $T_uM$. 
This implies that  $\tilde A_{\Sigma}$ is an abelian subspace of 
$\mathfrak {sim}^0(T_uM)$, the trace-free  symmetric endomorphisms of $T_uM$. 
In particular, $\tilde A_{\Sigma _j}$ is an abelian subspace of $\mathfrak {sim}^0(T_uM)$, for all $j=0, \cdots , r$. 
\

\begin {defi} The (restricted) normal holonomy factor $(\nu _i , \Phi _i)$ at $u$, $i>1$,  is said to have sections of the {\it compact type} if the abelian subspace $\tilde A_{\Sigma _i}$ of $\mathfrak {sim}^0 (T_uM$) is of  compact type, for any section $\Sigma _i \subset \nu _i$ of the action of $\Phi _i$ on $\nu _i$ (where $\tilde A$ is the traceless shape operator). 
\end {defi}

Observe that if $\Phi _i$ has sections of  the compact type, for some 
$u\in M$, then the same holds for any $u\in M$, since $M=H\cdot v$ is homogeneous. In this case we simply say that the normal holonomy factor $(\nu _i , \Phi _i)$ has sections of  the  compact type.

\begin {lema}\label {curvatura producto} Let $R$ be the curvature tensor at $p$ of a (symply connected) symmetric space $M$ without Euclidean factor. 
Assume that there are two complementary   subspaces  $\mathbb V, 
\mathbb W$ of  $T_pM$ such that $R_{\mathbb V, \mathbb W}=\{0\}$. 
Then  $\mathbb V \perp \mathbb W $ and so $M$ splits if the subspaces are non-trivial. 
\end {lema}

\begin {proof} From the Bianchi identity one obtains that 
\begin {equation} \label {1} 
R_{\mathbb V, \mathbb V}\mathbb W =\{0\}, \ 
R_{\mathbb W, \mathbb W}\mathbb V =\{0\}. 
\end {equation}
One has  that $R_{x,y}\cdot R =0$, where 
  $(R_{x,y}\cdot  R)_{u,v}= [R_{x,y},  R_{u,v}] 
  - R_{R_{x,y}u, v} - R_{u, R_{x,y}v}$. 
  So, from \eqref {1} one obtains that 
  \begin {equation} \label {2}
  [R_{\mathbb V, \mathbb V}, R_{\mathbb W, \mathbb W}] =0. 
  \end {equation}
  
  Assume that there exists  $z \in T_pM$ with 
  $z \notin \mathbb W$ and 
  $R_{\mathbb V, \mathbb V}
z=\{0\}$. Since $\mathbb V$ and $\mathbb W$ are complementary, we may assume, by \eqref {1}, that 
$z\in \mathbb V$. This implies, from the assumption $R_{\mathbb V, \mathbb W}=\{ 0 \}$ and \eqref {1},  that 
 $R _{T_pM, T_pM}z =\{0\}$. A contradiction since 
 	$R$ must  non-degenerate. The same is true if one interchanges $\mathbb V$ with $\mathbb W$.
 	
 	We have proved that 
 	
 	\begin {equation}\label {3}
 	\bigcap _{v, v' \in \mathbb V}\mathrm {ker}(R_{v, v' }) 
 	= \mathbb W,  \ \ \ 
 	\bigcap _{w, w' \in \mathbb W}\mathrm {ker}
 	(R_{w, w' }) 
 	= \mathbb V 
 	\end {equation}

 	Observe, from \eqref {2}, that any  skew-symmetric endomorphism $R_{v,v'}$,  $v,v' \in \mathbb V$,  commutes with all 
 	$R_{w,w'}$, 
 	$w,w'\in \mathbb W$. Then $R_{v,v'}$ leaves invariant 
 	$$\bigcap _{w, w' \in \mathbb W}\mathrm {ker}
 	(R_{w, w' }) $$
	 which coincides, by \eqref {3},  with $\mathbb V$ 
 	(and the same is true interchanging $\mathbb V$ with $\mathbb W$). Then 
 \begin {equation}\label {4}
 R_{\mathbb V , \mathbb V}\mathbb V \subset 
 \mathbb V, \  \ \ \ \ \ R_{\mathbb W , \mathbb W}\mathbb W \subset 
 \mathbb W 
 \end {equation}
 
 Observe that the image $\mathrm {im} (B)$ and the kernel  $\mathrm {ker}(B)$ of a skew-symmetric endomorphism $B$ are mutually perpendicular. Then, by making use of \eqref {4}, 
 the linear span of  $\{\mathrm {im} (R_{v,v'}): v,v'\in \mathbb V\}$ coincides with $\mathbb W ^\perp$. But, by 
 \eqref {4}, such a span is contained in $\mathbb V$. Then 
 $\mathbb W ^\perp \subset \mathbb V$.  Since
 $\dim (\mathbb W ^\perp) = \dim (\mathbb V)$, we conclude that $\mathbb V = \mathbb W ^\perp$. 
 
\end {proof}

\section {Proof of \normalfont {(b) $\Rightarrow$ (a) } of Theorem  \ref {main}}

Let $p\in M^n$ be fixed and, for simplicity of the exposition, let us identify 
$T_pM\simeq \mathbb R^n$ (by means of an orthonormal basis). Let us consider  the compact symmetric space $X^* = \mathrm {SU}_n/\mathrm {SO}_{n}$ with Cartan decomposition 
$\mathfrak{su} _n = i\, \mathfrak {sim}^0_n \oplus 
\, \mathfrak {so}_n$. Let us consider the map 
$ \tilde A^*: \nu _i \to i\, \mathfrak {sim}^0_n$, 
defined by 
$$\tilde A^*_\xi = i \tilde A_\xi, $$
where $\tilde A$ is the traceless shape operator of the submanifold $M$ of the sphere and $(\nu _i, \Phi _i)$ is an irreducible and non-transitive normal holonomy factor of the normal holonomy of $M$ at $p$.

From the Ricci identity, the adapted normal curvature tensor $\mathcal R$  (see \cite {BCO}) is given by
$$\langle \mathcal R _{\xi _1, \xi _2}\xi _3, \xi _4 \rangle = 
\textrm {trace}([\tilde A_{\xi _1}, \tilde A_{\xi _1}]   [\tilde A_{\xi _3}, \tilde A_{\xi _4}]),$$

 Then 
 $$\langle \mathcal R _{\xi _1, \xi _2}\xi _3, \xi _4 \rangle =-\langle 
  R_{\tilde A^*_{\xi _1}, \tilde A^*_{\xi _2}}\tilde A^*_{\xi _3}, \tilde A^*_{\xi _4} \rangle $$
  (see Remark \ref{curvature}).
  
  Let us consider the algebraic curvature tensor 
  $\bar R$ of $\mathbb V:= \tilde A^*_{\nu _i} \subset i\, \mathfrak {sim}^0_n $ given by 
  $$\langle \bar R _{u,v}w,z \rangle  = \langle  R _{u,v}w,z \rangle,  $$
  i.e., $\bar R _{u,v}w$ is the projection to $\mathbb V$ of $ R _{u,v}w$.
  
  Observe that $\bar R\neq 0$. Otherwise,  $\mathcal R _{\vert \nu _i} = 0$ and so $\nu _i$ would be flat (we identify the subspace $\nu _i $ of 
  $\nu _pM$ with the associated parallel sub-bundle of $\nu M$). 
  
  Since the factor $\Phi _i$ of $\Phi (p)$ acts irreducibly on $\nu _i$, and it  is non-transitive (on the sphere), 
  the map $\tilde A : \nu _i \to  \mathfrak {sim}^0_n$, and hence $\tilde A^* :
  \nu _i \to i\, \mathfrak {sim}^0_n$ is injective 
  \cite {OR}. In this reference it was used the shape operator $A$, but the proof is the same if one replaces $A$ by $\tilde A$. Moreover, one can  consider only $\tilde A$ restricted $\nu _i$.
  
  \begin {lema}\label {non-transitive}Let $\mathfrak h$ be the Lie subalgebra of $\mathfrak {so}(\mathbb V)$ which is algebraically generated by $\{ \bar R _{u,v}: u, v\in \mathbb V\}$ and let 
  $H\subset \mathrm {SO}(\mathbb V)$ be its associated Lie group. Then $[\mathbb{V}, \bar R, H]$ is an irreducible non-transitive holonomy system {\rm (see \cite {BCO}).} 
  \end {lema} 
  
  \begin {proof} Let $\mathcal O$ be the open and dense subset of $\nu _i$ which consists of the regular vectors for the  representation of $\Phi _i$ on $\nu _i \subset \nu _pM$. Let $\mathcal D$ be the (smooth) distribution on $\mathcal O$ given by the normal spaces to the $\Phi _i$-orbits (or, equivalently, to the sections of $\Phi _i$). Let $\bar {\mathcal O}: = \tilde A^*_{\mathcal O}$. Let us prove that the distribution $\bar {\mathcal D}^\perp$ is integrable. First, observe the following: the (isotropy) orbits of  $\text {SO}_n$ on 
  $i\, \mathfrak{sim}^0_n$ are always perpendicular to any abelian subspace. Then the integrability of $\bar {\mathcal D}^\perp$  follows from Proposition \ref {mainpro}.
  (Recall, from assumptions, that  $ A^* _{\Sigma}$ is an abelian subspace of compact type of any $\Sigma$ section of the normal holonomy action).
  
 Let $u,v\in \mathbb V$ and let us consider the (Euclidean) Killing field $X^{u,v}$ of $\mathbb V$ given by $X^{u,v}_q = \bar R _{u,v}q$. Then
 $X^{u,v}_{\vert \mathcal O}$ lies in $\mathcal D ^\perp$.  Since this distribution is integrable, any iterated bracket of Killing field of the form 
 $X^{u,v}_{\vert \mathcal O}$, $u,v \in \mathbb V$ lies in  $\mathcal D ^\perp$. 
 This implies that $H$ is non-transitive (on the unit sphere), since $\dim \mathcal D\geq 2$.
 
 Since $[\nu_ i; \mathcal R, \Phi _i]$ is an  irreducible holonomy system and 
 $$\langle \mathcal R _{\xi _1, \xi _2}\xi _3, \xi _4 \rangle =-\langle 
 \bar R_{\tilde A^*_{\xi _1}, \tilde A^*_{\xi _2}}\tilde A^*_{\xi _3}, \tilde A^*_{\xi _4} \rangle, 
 $$ 
 one obtains that if the holonomy system $[\mathbb{V}, \bar R, H]$ is reducible, there would exist a non-trivial orthogonal decomposition 
$\mathbb V = \mathbb V _1\oplus  \mathbb V _2$ such that 
$\bar R _{\mathbb V _1 , \mathbb V_2} = \{0\}$.
Then the subspaces $(\tilde A^*)^{-1}(\mathbb V _1)$ and $(\tilde A^*)^{-1}(\mathbb V _2)$ are in the assumptions of Lemma \ref {curvatura producto} and $\mathcal R$ is a non-trivial product of Riemannian curvature tensors. A contradiction. 
This shows the irreducibility of $[\mathbb{V}, \bar R, H]$.

\end {proof}

By Lemma \ref {non-transitive} and \cite [Proposition 2.21] {OR} we obtain 
the following:

\begin {cor}\label{coroholo} $\tilde A^*: \nu _i \to \mathbb V = \tilde A^*_{ \nu _i}$ is a homothety and $ \tilde A^* \Phi _i  (\tilde A^*)^{-1} = H$. 

\end{cor}

From the above corollary and Proposition \ref {mainpro} one obtains, by a limit argument for singular points,  that: 

\

\noindent {\it {\bf (*)}  given $\xi \in \nu _i$, $\phi \in \Phi _i$ there exists 
$k\in K'= \mathrm {SO}_n$ such that 
$\tilde A^*_{\phi \xi} = k \tilde A^*_{\phi \xi}k^{-1} $}

\

(recall that $\mathrm {SO}_n$ acts on the symmetric matrices by conjugation)

\

Then $\tilde A^*_{\phi \xi}$ and $\tilde A^*_{\xi}$ have the same (pure imaginary) eigenvalues. The same is true if we replace $\tilde A^*$ by the shape operator $A$. 
Namely, $\tilde A^*_{\phi \xi}$ and $\tilde A^*_{\xi}$ have the same eigenvalues. 

Being the homogeneous submanifold $G\cdot p =M^n$ an irreducible submanifold of $\mathbb R^N$, by \cite [Theorem 5.2.4] {BCO}, given $g\in G$ there exists a curve 
$c$ in $M$ from $p$ to $gp$ such that $\mathrm {d}g_{\vert \nu _i}$ coincides with the normal parallel transport $\tau  _c$ along $c$. This property, together with (*), implies that for any curve $\gamma $ in $M$ and any 
$\xi \in (\nu _i) _{\gamma (0)}$ $\tilde A_\xi$ and $\tilde A_{ \tau  _{\gamma}(\xi )}$ have the same eigenvalues.  Then, by Theorem \ref {traceless} and Remark \ref {subbundle},  $M$ has constant principal curvatures. Since the normal holonomy of $M$ has a non-transitive normal holonomy factor,  any principal  holonomy tube of $M$ has codimension at least $2$ in the sphere. Then, by the Theorem of Thorbergsson \cite {Th} (see also {BCO}), one has that $M$ is an orbit of an irreducible $s$-representation. This completes the proof of 
 (b) $\Rightarrow$ (a) of the Theorem \ref {main}.

\section {The proof of \normalfont {(a) $\Rightarrow$ (b)} of Theorem  \ref {main}}

\

Let $ G/K$ be an irreducible simply connected symmetric space with 
Cartan decomposition $\mathfrak g = \mathfrak k \oplus \mathfrak p$. 
Let us consider the isotropy representation of $K$ on $T_{[e]}(G/K)$ which can be regarded as 
 the $\mathrm {Ad}$-representation of $K$ on $\mathfrak p\simeq \mathbb R^n$, whose scalar product will be denoted by $(\, ,\, )$. 
 We will regard $\mathrm {Ad}(K)$ as a compact subgroup, that we will also denote by $K$,  of 
 $\mathrm {SO}(n)$.  Let us consider a principal orbit $K\cdot p$ which is an isoparametric submanifold. Let $\Sigma$ be the normal space at $p$ of 
 $K\cdot p$ that corresponds to a maximal abelian subspace of the Cartan complement $\mathfrak p$. Let $W$ be the (irreducible)  Weyl group at $p$ associated to the isoparametric submanifold $K\cdot p$. Note that $W$ does not change if we pass to a parallel principal orbit $K\cdot q$, $q\in \Sigma$. The finite reflection group $W$ coincides with the usual one associated to the symmetric space $G/K$, with respect to the maximal abelian subspace $\Sigma$. Such a Weyl group $W$   has associated a reduced  (crystallographic) root  system $\Phi \subset \Sigma$, whose elements are called roots. That is, $\Phi$ is a finite subset of $\Sigma$ such that the following conditions holds: 
 \begin {enumerate}
 \item $0\notin \Phi$ and $\Phi$ spans $\Sigma$.
 
 \item If $x \in \Phi$ then $-x\in \Phi$ and no other scalar multiple of 
 $x$ belongs to $\Phi$.
 
 \item The reflection through  the hyperplane perpendicular to any root leaves 
 $\Phi$ invariant.

 \item The number $2\frac {(x,y)}{(x,x)}$ is a rational number, for any $x, y\in \Phi$.
 
 \end {enumerate}
 
 \
 
 The Weyl group $W$ coincides with the finite group of isometries which is generated by the reflections through the hyperplanes perpendicular to the roots. 
 
 Observe that from (iv) we deduce that $\frac {(y,y)}{(x,x)}$ is a rational number, provided $(x,y)\neq 0$. But this will be always true if 
 $W$ acts irreducibly. In fact, given $x,y\in \Phi$ there exists a  finite sequence of elements $x_1, \cdots , x_j$ in $\Phi$ such that 
 $x_1 = x$, $x_j=y$ and such that $(x_i,x_{i+1})\neq 0$, for $i=1, \cdots , j-1$. This implies that always $\frac {(y,y)}{(x,x)}$ is a rational number. (Otherwise the root system would be the union of two mutually orthogonal subsets and $W$ would be reducible). 
 
 Then, by rescaling all the roots by an appropriate $\lambda \neq 0$, we 
 may assume that the square of the norms of any element of $\Phi$ is a rational number. Then,  by (iv), the scalar product of any two elements of 
 $\Phi$ is a rational number.  
 
 Let us consider the discrete  abelian  subgroup  $L\subset \Sigma$ given by the linear combinations with integer coefficients of the  roots  
 (which is a so-called rational lattice). Let us consider a subgroup  $S$ of $L$ and  assume that 
 $S$ (linearly) spans a proper  subspace 
 $\mathbb W$ of $\Sigma$. Consider 
 $S^\perp = L\cap \mathbb W ^\perp$. Then 
  $S^\perp$ spans $\mathbb W^\perp$, i.e., it is a full lattice of $\mathbb W^\perp$. In fact, let $v_1 , \cdots ,  v_j \in S$ linearly independent that  spans $\mathbb W$ and complete them, with elements of $L$,    $v_{j+1}, \cdots ,  
  v_k$  to a basis of $\Sigma$. Use Gram-Schmidt procedure, but without normalizing the obtained vectors to obtain an orthogonal (no orthonormal) basis 
  $\tilde v_1, \cdots , \tilde v_k$ of $\Sigma$ such that any element of 
  $\{\tilde v_1, \cdots , \tilde v_k\} $ is a rational combination of elements of $L$ and 
  the span of $\tilde v_1, \cdots , \tilde v_i$ coincides with that of $ v_1, \cdots  v_i$,  $i=1, \cdots k$.
  We have used that the scalar product of any two elements of $L$ is rational.  After multiplying each of  the elements of    $\{\tilde v_1, \cdots , \tilde v_k \} $, by an appropriate integer number,  we obtain an orthogonal basis   $v'_1, \cdots \tilde v'_k $ of $\Sigma$ which consists of elements of $L$. It is clear that 
  $v'_{j+1}, \cdots , v'_k$ belong to $S^\perp$ and span $\mathbb W ^\perp$. 
 This may be interpreted as follows:  the exponential of the normal space at a point of  any subtorus of $\Sigma/L$ is a subtorus (cf. \cite {HPTT}[Thm. 2.8]).
 
 As a consequence, we obtain the following: {\it  if a  (positive dimensional) subspace $\mathbb V$ of $\Sigma$ is the intersection of some 
 of  the  reflection hyperplanes, then the elements of the lattice that belong to $\mathbb V$ span $\mathbb V$.}   In fact, let  $J\subset \Phi$ be  the roots associated to these hyperplanes.  Then $\mathbb V$ is the orthogonal complement of the linear span of $J$. 

\

 \subsection {The shape operator of isotropy orbits.}

\

\vspace {.2cm}

 We keep the assumptions and notation of this section.  Let $K\cdot p\subset \mathbb R ^N$ be a principal orbit, where $p\in \Sigma$. For $\alpha \in \Phi$ let us denote by 
 $H_\alpha$ the hyperplane of $\Sigma$ which is perpendicular to $\alpha$. 
 The set of  curvature normals $\{\eta _\alpha :\alpha \in \Phi\} $ at $p$, associated the commuting family of shape operators of $K\cdot p$, and indexed by $\Phi$,  are given by (see \cite {PT}): 
$$ \eta _\alpha = -\frac {1}{(p, \alpha)}\alpha .$$
So the eigenvalues of the  shape operator 
$A_{\xi }$, if $\xi \in \nu _p(K\cdot p) \simeq \Sigma$,  are given by 
 $-\frac {(\alpha ,\xi)}{(p, \alpha)}$. Let $p\in L$ be fixed. Then, if 
 $\xi \in \Phi$, the eigenvalues of 
 $A_\xi$ are rational numbers. Since $\Phi$ spans $\Sigma$, we obtain that the family of 
 shape operators of $K\cdot p$ is an abelian subspace of compact type of $\mathfrak {sim}(T_p(K\cdot p))$  (see Definition \ref {compact-type}). This family of shape operators is not in general of compact type if $p$ does not belong to the lattice $L$.
 
 Let us now assume that $K\cdot p$, $p\neq 0$, is not a principal orbit and so $\Sigma$ is properly contained in the normal space 
 $\nu _p(K\cdot p)$.  Let 
 $$J:= \{\beta \in \Phi: p \in H_\beta\}\subset \Phi. $$ Then 
 $p$ belongs to the intersection of the hyperplanes $H_\beta$, with $\beta \in J$. 
 Such an intersection is the orthogonal complement of  $\mathbb W := \text {linear span  of }J$. Assume, furthermore, that $K\cdot p$ is most singular. Namely, that 
 $\mathbb W^\perp =\mathbb R p$ is one dimensional.
 
Observe that $\mathbb W$ is  a section of the normal holonomy of  $K\cdot p$, regarded as a submanifold of the sphere.

 Let us compute the eigenvalues of the commuting family of shape operators 
 $\{A_\xi :\xi \in \Sigma\}$ of 
 $K\cdot p$ at $p$. One has that these eigenvalues are given by 
 
$$ \big\{-\frac {(\alpha ,\xi)}{(p, \alpha)}: \alpha \in \Phi \text { is not perpendicular to } p \big\}$$

 As we have remarked $\mathbb W^\perp$ has associated a full sub-lattice. So there must  exist a non-zero $\gamma  \in L$ such that $\gamma = a p$, for some  $a\in \mathbb R$.

Let us consider the set $\tilde J := a J $ which spans $\mathbb W$. 
Then any shape operator of $A_\xi$, with $\xi \in \tilde J$ has rational eigenvalues. So, $\tilde A_\xi$ has rational eigenvalues.
Then $\tilde A_\mathbb W \subset \mathfrak {sim}^0(T_p(K\cdot p))$ is of compact type. 
Since the isotropy $K_p$, represented on the normal space $\nu _p(K\cdot p)$, coincides with the normal holonomy, we obtain that the normal holonomy of $K\cdot p$ has sections of compact type. 

Let $\Phi$ be the (restricted) normal holonomy group of $K\cdot p$ at $p$ (regarded as a submanifold of the sphere). Let  us decompose, as in the normal holonomy theorem, 
$\Phi = \Phi _1 \times \cdots \times \Phi _r$,  $\nu _p (K\cdot p)= \nu _1 \oplus \cdots \oplus  \nu _r$, where $\Phi _i$ acts irreducible on 
$\nu _i$ and trivially on $\nu _j$, if $i\neq j$. Choose any normal holonomy factor that we may assume that it is 
$(\nu _1 , \Phi _1)$. Let $\Sigma$ be a section for the normal holonomy action of $\Phi$ on $\nu _p(K\cdot p)$ and let 
$\Phi ^1 = \prod _{i\neq 1}\Phi _i$. Then 
$$\Sigma ^1 = \bigcap _{g\in \Phi ^1}\, g\Sigma$$ 
is a section for the action of $\Phi _1$ on $\nu _1$. Since $\xi \mapsto \tilde A _\xi$ is injective 
$$\tilde A _{\Sigma ^1} =  \bigcap _{g\in \Phi ^1}\tilde A_{g\Sigma}= \bigcap _{g\in \Phi ^1}g\tilde A_{\Sigma}g^{-1}$$. 

This shows that $\Sigma ^1$ is a section of $\Phi _1$ of compact type. Hence any normal holonomy factor  $( \nu _1 , \Phi  _1)$, of 
the normal holonomy of $K\cdot p$,  has sections  of the compact type. This finishes the proof of $(a)$ implies $(b)$ of Theorem \ref {main}.

\qed
 
\

We do not know any homogeneous example that shows that the non-transitivity of the normal holonomy could not be removed from Theorem \ref {main}. But we expect it does exist. 

\begin {nota}\label {subbundle} Let $M^n$ be  an irreducible and full submanifold of the Euclidean space. Assume that the (restricted) normal holonomy $\Phi (p)$ of $M$ has a 
(non-necessarily irreducible) non-transitive factor. Namely, there is a subspace $\nu '_pM \subset \nu _pM$ which is 
$\Phi (p)$-invariant and such that $\Phi (p)$ is not transitive on the unit sphere of   $\nu '_pM$.
This defines a $\nabla ^\perp$-parallel sub-bundle $\nu 'M$ of the normal bundle $\nu M$.   Assume  that the shape operator $A_{\xi (t)}$ has constant eigenvalues for any arbitrary   parallel normal field $\xi (t)$ (along any curve in $M$ ) which lies in $\nu 'M$. 
Then $M$ has constant principal curvatures. In fact, if $w \neq 0$ belongs to 
$\nu '_pM$,  the holonomy tube $(M)_w$ admits  a parallel normal field 
$\tilde \eta$ such that $M= ((M)_w)_\eta $ (i.e., $M$ is a parallel focal manifold to the holonomy tube). From the tube formula (see \cite {BCO}),  the eigenvalues 
of the shape operator $A_\eta$ of $(M)_w$ are constant (and there are at least two different eigenvalues). Then the irreducible submanifold $(M)_v$ admits a so-called parallel isoparametric section. Then, by theorems 4.5.10  and  4.5.2 of \cite {BCO}, $(M)_w$, and hence $M$  is contained in a sphere and has constant principal curvatures.

\

\section {The proof of Theorem \ref {traceless}}

This proof is independent of the main tools developed in this article. 

\

\begin {proof} If $M$ is a submanifold of constant principal curvatures, then 
its  mean curvature normal field is parallel and so the traceless shape operator has constant principal curvatures, and in particular when 
restricted to $\nu'M$. Let us prove the converse.

If $M$ is a submanifold of the sphere $S^{N-1}$, we will regard it as an Euclidean submanifold. In this case its normal space $\nu M$ will be regarded as the parallel sub-bundle, of  the normal space to the $\mathbb R^N$, which is perpendicular to the position vector field. 

Let $\Sigma \subset \nu'_pM$ be section for the normal holonomy action of $\Phi (p)_{\vert \nu'_pM}$ at $p\in M$. Then the family of traceless shape operators $\tilde A_\Sigma = \{\tilde A_\xi : \xi \in \Sigma\}$ 
is simultaneously diagonalizable. Each eigenvalue linear function $\lambda ^\Sigma : \Sigma \to \mathbb R$  defines the so-called associated 
curvature normal  $\eta ^\Sigma $ (with respect to $\nu'M$ and $\Sigma$)
, i.e.,  $\langle \eta^\Sigma, \, \cdot \,\rangle = \lambda (\,\cdot \,)$. Let $\eta ^\Sigma _1, \cdots , \eta ^\Sigma _d$ be the 
distinct curvature normals. Let $c(t)$ be a curve on $M$ with $c(0) =p$ and let $\tau _t$ denote the $\nabla ^\perp$-parallel transport along $c(t)$, restricted to $\nu' M$. 
If $\xi (t)$ is a parallel normal field along $c(t)$, with $\xi (0)\in \Sigma$,   then $\eta^{\tau _t(\Sigma)}_ i = \tau _t(\eta ^\Sigma _i)$ are the curvature normals associated to 
$\tilde A_{\tau _t(\Sigma)}$, $i=1, \cdots ,  d$.
Let  $i _0\in \{1, \cdots , d\}$ be fixed. Let us show that 
$\eta_i^\Sigma - \eta_{i_0}^\Sigma$, $i=1, \cdots , d$, span
$\Sigma$. If not, let $0\neq \xi \in \Sigma$ perpendicular to this span. This implies that $\tilde A_\xi =0$ and so $A_\xi = \langle \frac 1 n H(p), \xi \rangle 
\mathrm {Id} =\langle \frac 1 n H'(p), \xi \rangle 
\mathrm {Id}$ , where $H'$ is the projection of $H$ to the parallel subbundle $\nu'M$. From the assumptions, if $\xi (t)$ is the parallel   transport  in the normal connection along any curve $c(t)$, then $\tilde A _{\xi (t)} = 0$. 
So $A_{\xi (t)} =  \langle \frac 1 n H'(c(t)), \xi (t) \rangle 
\mathrm {Id}$. Let $\bar \nu M$ be the parallel subbundle of the parallel  subbundle $\nu' M\subset \nu M$ that is linearly spanned by the normal parallel transport of $\xi $ along any curve. Then $A_\psi = \langle \frac 1 n H'(q), \psi \rangle 
\mathrm {Id}$, for any $\psi \in \bar \nu' _qM$. The same proof, relaying on the Codazzi identity, used for proving that umbilical submanifolds have parallel mean curvature, shows that the orthogonal projection $\bar H$ of $H'$ to $\bar \nu M$ is parallel. Assume that $M$ is a Euclidean submanifold. Then, since $A_{\bar H}$ is a constant multiple of the identity, 
$M$ is either contained in a sphere, if $\bar H\neq 0$, or $M$ is not full. This contradicts the assumptions. 
If $M$ is a spherical submanifold, then in any case $M$ is not full. A contradiction. Then $\eta_i^\Sigma - \eta_{i_0}^\Sigma$, $i=1, \cdots , d$, span
$\Sigma$.

Let  $p\in M$ be arbitrary  and let  $\Sigma $ be a  section for  $\Phi (p)_{\vert \nu'_pM}$.  Let $i_0 \in \{1, \cdots , d\}$ be fixed and let, for $i \in \{1, \cdots , d\}$,  $\mathbb V_i = \{\eta_i^\Sigma - \eta_{i_0}^\Sigma\}^\perp$. 
If the hyperplanes $\mathbb V_i $, $i\neq i_0$, do not generate $\Sigma $ then 
$\eta_i^\Sigma - \eta_{i_0}^\Sigma$ must all lie in a line through the origin of 
$\Sigma $. This is a contradiction since $\eta_i^\Sigma - \eta_{i_0}^\Sigma$ span $\Sigma $ and $\dim \Sigma \geq 2$. 

Let $v_i\in \mathbb V_i$ be fixed and short and let us consider the (partial) 
holonomy tube \cite [3.4.3] {BCO}
$$(M)_{v_i} =\{c(1)  + \tau _c  (v_i)\}$$
where $c:[0,1]\to M$ with $c(0)=p$ and $\tau _c $  denotes  the normal parallel transport
  along  $c$ (eventually, by making $M$ smaller around $p$).
 
 We have that $M$ itself is parallel focal manifold to the holonomy tube. 
 Namely, 
 $$M= ((M)_{v_i})_{-\xi_i}$$ 
 where $\xi _i$ is the parallel normal field of 
 $(M)_{v_i}$ defined by $\xi _i(c(1) + \tau _c(v_i)) = \tau _c(v_i)$. 
 Let us consider the shape operator $A^i_{\xi _i}$ of $(M)_{v_i}$. 
 Since $M= ((M)_{v_i})_{-\xi_i}$, then the vertical distribution 
 $$\nu ^i := \ker (\mathrm {Id} - A^i_{-\xi _i})$$
  of $(M)_{v_i}$
 coincides with the eigendistribution of 
  $A^i_{-\xi _i}$ associated to the eigenvalue $1$. The restriction of
   $A^i_{\xi _i}= - A^i_{-\xi _i}$ to the horizontal distribution $\mathcal H$, perpendicular to $\nu ^i$, is given by the so-called tube formula  in \cite {BCO}. 
   
   Namely,  let  $q= c(1) +  \tau _c(v_i) $. Then, taking into account that $\mathcal H _q = T_{c(1)}M$, 
   
   $$A^i_{\xi _i (q) \vert T_{c(1)}M}= 
   A_{\tau _c(v_i)}(Id-  A_{\tau _c(v_i)})^{-1}.$$
and no eigenvalue of this restriction is $-1$. 

From the assumptions,  the eigenvalues of the traceless shape operator 
$\tilde A _{\tau _c(v_i)}$ are the same as those of $\tilde A _{v_i}$. Then  the multiplicities of the eigenvalues 
of $A_{\tau _c(v_i)}$ do not depend on $c$. This implies that the multiplicities of the eigenvalues of $A^i_{\xi _i}$ are constant. Let $\mathbb W _1, \cdots , \mathbb W _d$ be the common eigenspaces of $A_{\Sigma }$,  associated to the eigenvalues functions $\langle \eta _1^{\Sigma } , \, \cdot \, \rangle, \cdots , \langle \eta _d^{\Sigma } , \, \cdot \, \rangle$. 
   
Since $v_i$ is perpendicular to $\eta _i^{\Sigma }-\eta _{i_0}^{\Sigma }$, 
$\langle \eta _i^{\Sigma }, v_i\rangle =\langle \eta _{i_0}^{\Sigma }, v_i\rangle$.

 From the tube formula, we obtain that  $\mathbb W_{i_0} \oplus \mathbb W_i$ is contained in an eigenspace of 
$A^i_{\xi _i(q)}$ which extends to an (integrable by Codazzi identity) eigendistribution $E^i$ of $A^i_{\xi _i}$, since the eigenvalues of this shape operator have constant multiplicities. Moreover, $E^i$ is perpendicular to the vertical distribution $\nu ^i$.
Since $\dim (E^i)\geq 2$, we have that the corresponding eigenvalue, let us say $\lambda _i$, is constant along $E^i$ (this is similar to the classic Dupin condition; see \cite [Lemma 3.3]{OR}). This eigenvalue is given by 

$$\lambda _ i (q) = \frac {b_i 
+ \frac 1 n \langle H'(c(1))  , \tau _c(v_i)\rangle }
{1- b_i
	- \frac 1 n \langle H'(c(1))  , \tau _c(v_i)\rangle}\ , $$
where 
$$b _i : = \langle \eta ^{\Sigma }_{i_0}(p), v_i \rangle  =  \langle \eta ^{\tau _c (\Sigma )}_{i_0}(c(1)), \tau _c(v_i)\rangle. $$ 

Now let $\tilde c (t)$ be a curve in $(M)_{v_i}$ which is tangent to $E^i$ with $\tilde c (0) = p+ v_i$. Let $c(t) = \pi (\tilde c (t))$, where 
$\pi: (M)_{v_i} \to M$ is the natural projection 
$u\overset {\pi }{\to}u -\xi _i(u)$. Then $v_i(t) = \tilde c (t) -c(t)$ coincides with the normal parallel transport in $M$ of $v_i$ along $c(t)$.  Then the previous formula yields 

$$\lambda _ i (\tilde c (t)) = \frac {b_i 
	+ \frac 1 n \langle H'(c(t))  , v_i(t)\rangle }
{1- b_i
	- \frac 1 n \langle H'(c(t))  , v_i (t)\rangle}\ . $$

The constancy of $\lambda _i$ along $E^i$ implies, by taking the derivative,  that $\nabla ^\perp _{c'(0)}H'$ is perpendicular to $v_i$.
Observe now that 

$$c'(0) = \frac {\, \mathrm d}{\mathrm dt}_{\vert 0} (\tilde c (t)  -\xi _i(\tilde c (t)) = ( Id+A^i_{\xi _i (p+v_i)})\tilde c' (0) = 
( Id+A^i_{v_i})\tilde c' (0) .$$

Observe that $A^i_{v_i}E^i(p+v_i)=E^i(p+v_i)$. So, since $\tilde c(t)$ is arbitrary, we obtain that $\nabla ^\perp _{E^i(p+v_i)}H'$ is perpendicular to $v_i$. Hence, since $\mathbb W_{i_0}\subset \mathbb W_{i_0}\oplus \mathbb W_{i}\subset 
E^i(p+v_i)$, 

$$ \displaystyle { v_i \perp \nabla ^\perp _{\mathbb W_{i_0}}H' }\ .$$

 Since $i\neq i_0$ is arbitrary, and  $v_i$ is arbitrary in $\mathbb  V_i$, and these hyperplanes generate $\Sigma$, we obtain that  $\nabla ^\perp _{\mathbb W_{i_0}}H'\perp \Sigma$. Since $i_0$ is arbitrary, we have that 
$\nabla ^\perp _{\mathbb W_{i}}H'\perp \Sigma$, for all $i=1, \cdots , 
d$. Then $\nabla ^\perp _{T_pM}H'\perp \Sigma$. Since $\Sigma$ is arbitrary section of $\Phi (p)_{\vert \nu'M}$, and since any vector 
in $\nu' _p M$ is contained in some section, we conclude that 
$\nabla ^\perp _{T_pM}H' =0$. Since $p\in M$ is arbitrary, we conclude that 
 $H'$ is parallel. This implies, by making use the assumptions,  that the shape operator of $M$, restricted to  $\nu'M$, has constant principal curvatures along parallel normal fields. Thus, by Remark \ref {subbundle}, $M$ has constant principal curvatures.

\end {proof}

\end {nota}

\end {document}